\documentclass[11pt]{amsart}
\usepackage[T1]{fontenc}
\usepackage{lmodern,microtype,amsmath,amssymb,amsthm,mathtools,mathrsfs,booktabs}
\usepackage[margin=27mm]{geometry}
\usepackage[hidelinks]{hyperref}
\newtheorem{theorem}{Theorem}[section]
\newtheorem{proposition}[theorem]{Proposition}
\newtheorem{lemma}[theorem]{Lemma}
\newtheorem{corollary}[theorem]{Corollary}
\theoremstyle{remark}
\newtheorem{remark}[theorem]{Remark}
\newcommand{\R}{\mathbb R}
\newcommand{\D}{\mathcal D}
\newcommand{\Div}{\operatorname{div}}
\newcommand{\tr}{\operatorname{tr}}
\newcommand{\supp}{\operatorname{supp}}
\newcommand{\capc}{\operatorname{Cap}_{2,2}}
\newcommand{\lamD}{\lambda_1^D}
\newcommand{\calR}{\mathscr R}
\newcommand{\calA}{\mathscr A}
\newcommand{\calB}{\mathscr B}
\newcommand{\calD}{\mathscr D}
\newcommand{\calE}{\mathscr E}
\newcommand{\Q}{\mathcal Q}
\numberwithin{equation}{section}
\allowdisplaybreaks[1]
\title{Rayleigh's conjecture in dimensions four and five}

\author[D. Suragan]{Durvudkhan Suragan}
\address{
	Durvudkhan Suragan:
	\endgraf
	Department of Mathematics
	\endgraf
	Nazarbayev University
	\endgraf
	Astana 010000
	\endgraf
	Kazakhstan
	\endgraf
	{\it E-mail address:} {\rm  durvudkhan.suragan@nu.edu.kz}
}

\begin{document}
\maketitle

\begin{abstract}
We prove Rayleigh's conjecture for the clamped plate in dimensions four
and five. In fact, we obtain a sharp comparison with the first Dirichlet
eigenvalue: on every open set of finite positive measure, the first plate
eigenvalue is at least the square of the Dirichlet eigenvalue multiplied
by the ratio attained on a ball. The proof uses a homogeneous
transformation of a trial function and its gradient. Two divergence
identities reduce the comparison to positivity of a quadratic form in the
Hessian. We establish the required Bessel inequalities by rational
estimates for entire profiles and their remainders. An
exact defect identity accounts for regular nodal sets, so the argument
allows trial functions of either sign and imposes no boundary regularity.
The domains may be disconnected or unbounded.
We characterize equality up to sets of zero Sobolev capacity. We also
obtain quantitative stability with the optimal quadratic power of the
Fraenkel asymmetry.
\end{abstract}

\noindent\textit{2020 Mathematics Subject Classification.}
35P15, 35J40, 49Q20.\\
\textit{Keywords.} Clamped plate, Rayleigh conjecture, spectral comparison,
Bessel functions, quantitative isoperimetry.

\section{Introduction and the main result}

Let $\Omega\subset\R^n$ be open with $0<|\Omega|<\infty$, where
$|E|$ denotes Lebesgue measure. For $k=1,2$, let $H^k(\R^n)$ be the
usual Sobolev space of functions whose weak derivatives through order
$k$ belong to $L^2(\R^n)$. We define $H_0^k(\Omega)$ as the closure
in $H^k(\R^n)$ of the smooth, compactly supported functions
$C_c^\infty(\Omega)$, using zero extension. All trial functions are
real-valued.
The first eigenvalues of the clamped biharmonic operator and the
Dirichlet Laplacian are
\begin{align}
 \Gamma_1(\Omega)&=\inf_{0\ne U\in H_0^2(\Omega)}
 \frac{\int_{\R^n}(\Delta U)^2}{\int_{\R^n}U^2},\label{eq:Gamma}\\
 \lamD(\Omega)&=\inf_{0\ne v\in H_0^1(\Omega)}
 \frac{\int_{\R^n}|\nabla v|^2}{\int_{\R^n}v^2}.
 \label{eq:Dirichlet}
\end{align}
Here $\nabla$ and $\Delta$ denote the gradient and the Euclidean
Laplacian. On a smooth bounded domain, the first problem has the boundary
conditions $U=\partial_{\mathbf n}U=0$, where $\partial_{\mathbf n}$
is the outward normal derivative. Rayleigh's conjecture asserts that
the ball minimizes $\Gamma_1$ at fixed volume. We prove the conjecture
in dimensions four and five by comparing the two eigenvalues on the
same domain.

Let $J_\sigma$ and $I_\sigma$ be the ordinary and modified Bessel
functions of the first kind. Write $j_{\sigma,k}$ for the $k$-th
positive zero of $J_\sigma$. Put $\nu=n/2-1$ and $q_n=j_{\nu,1}$,
and let $\beta_n$ be the first positive zero of
\begin{equation}\label{eq:crossproduct}
 J_\nu(r)I_{\nu+1}(r)+J_{\nu+1}(r)I_\nu(r).
\end{equation}
Write $B_R(a)=\{x\in\R^n:|x-a|<R\}$ and $B_R=B_R(0)$.
Set $\omega_n=|B_1|$ and
$\Omega^*=B_{(|\Omega|/\omega_n)^{1/n}}$, the centered ball of
volume $|\Omega|$.

\begin{samepage}
\begin{theorem}\label{thm:main}
Let $n\in\{4,5\}$. Every open set $\Omega\subset\R^n$ with
$0<|\Omega|<\infty$ satisfies
\begin{equation}\label{eq:comparison}
 \Gamma_1(\Omega)\ge
 \left(\frac{\beta_n}{q_n}\right)^4\lamD(\Omega)^2.
\end{equation}
The constant is sharp and is attained on balls. Consequently,
\begin{equation}\label{eq:Rayleigh}
 \Gamma_1(\Omega)\ge
 \beta_n^4\left(\frac{\omega_n}{|\Omega|}\right)^{4/n}
 =\Gamma_1(\Omega^*).
\end{equation}
Equality in \eqref{eq:Rayleigh} holds if and only if, after a translation,
$\Omega=B_R\setminus E$, where $E$ is relatively closed in $B_R$ and
$\capc(E)=0$. In particular, among bounded domains with $C^1$ boundary,
equality holds if and only if $\Omega$ is a ball.
\end{theorem}
\end{samepage}

The capacity in the theorem is the Sobolev capacity associated with
$H^2(\R^n)$; its definition and role in the equality case are given in
Section~\ref{sec:equality}. Table~\ref{tab:constants} gives the
relevant constants. These decimals are for orientation;
Section~\ref{sec:scalar} gives the rational enclosures used in the proof.
\begin{table}[htbp]
\centering
\begin{tabular}{cccc}
\toprule
$n$&$q_n$&$\beta_n$&$(\beta_n/q_n)^4$\\
\midrule
$4$&$3.83\ldots$&$4.61\ldots$&$2.09\ldots$\\
$5$&$4.49\ldots$&$5.26\ldots$&$1.88\ldots$\\
\bottomrule
\end{tabular}
\caption{Ball eigenvalue constants in dimensions four and five.}
\label{tab:constants}
\end{table}

An obstacle to a direct rearrangement argument is that a first
clamped-plate eigenfunction can change sign. Moreover, Schwarz
symmetrization need not decrease the norm of its Laplacian. Szeg\H{o}
\cite{S50,S58} treated the fixed-sign case. Talenti's comparison theorem
\cite{T76,T81} led to a relaxation on two coupled balls. Nadirashvili
\cite{N95} proved the planar conjecture, and Ashbaugh and Benguria
\cite{AB95} established the result in dimensions two and three.
Ashbaugh and Laugesen \cite{AL96} showed why this relaxation cannot prove
the higher-dimensional conjecture: for $n\ge4$ its minimum occurs at
two equal balls. Their lower bound is
\[
 \Gamma_1(\Omega)>d_n\Gamma_1(\Omega^*),\quad
 d_n=2^{4/n}(q_n/\beta_n)^4<1.
\]
In dimension four, $d_4=0.95\ldots$. Surveys of the plate and
buckling problems include \cite{ABL97,BBL}. Positivity for higher-order
operators is treated in \cite{GGS10}.

Several recent results address different parts of this difficulty.
Leylekian \cite{L24critical} proved ball optimality for sufficiently
regular minimizers in dimensions four through nine. His result assumes
a condition on the critical values of the first eigenfunction. Other
sufficient conditions involve its mean or an additional boundary
condition \cite{L24sufficient}. An open minimizing shape exists in
dimensions two through eight \cite{L25existence}. Its boundary
regularity remains an obstacle to applying those conditional results.
Ashbaugh, Bucur, Laugesen, and Leylekian \cite{ABLL25} proved
fourth-order Saint-Venant inequalities. In the preprint
\href{https://arxiv.org/abs/2610.29840}{arXiv:2610.29840}, Bucur and
Laugesen give a simpler planar proof and nonlinear extensions. Their
higher-dimensional theorem concerns an $L^\infty$ displacement norm and
leaves the classical $L^2$ problem open for $n\ge4$. Frank, Laptev,
and the author \cite{FLS26} give geometric lower bounds in terms of a
generalized inradius.

Our proof uses the radial plate and membrane modes of a ball to construct
a single-variable function $f$. For a real clamped trial function $U$, the
associated Dirichlet trial function is
\begin{equation}\label{eq:transformintro}
 F[U]=|\nabla U|f\left(\frac{|U|}{|\nabla U|}\right),
\end{equation}
with a continuous extension at critical points. The homogeneous
transformation and calibration strategy are similar to a construction recently appeared in 
\href{https://clopen.solutions/submissions/proof-of-the-planar-polyaszego-buckling-conjecture-2/v1/}{\emph{Proof of the planar P\'olya--Szeg\H{o} buckling conjecture}}
(unpublished preprint, 2026) by Tang
and Zhang. However, the plate energy requires a different
potential correction. In dimensions greater than two, a tangential
trace-free Hessian term also appears.

Two divergence terms convert the difference of the plate and membrane
energies into a quadratic form in $D^2U$. Its positivity is a question
about explicit radial profiles. The coordinate
$x=1-r^2/\beta_n^2$ expresses the coefficients through entire functions
and removes their endpoint cancellations. Rational Taylor bounds, including
the full tails, then give positivity throughout the interval. All estimates
needed for this step appear in Section~\ref{sec:global}.

The argument does not discard nodal sets. Proposition~\ref{prop:defect}
gives an exact identity whose nodal contribution is a positive multiple
of $\int_{\{U=0,\,|\nabla U|>0\}}|\nabla U|^2\,d\mathcal H^{n-1}$,
where $\mathcal H^{n-1}$ is $(n-1)$-dimensional Hausdorff measure.
This allows arbitrary real trials and requires no boundary regularity.

The Fraenkel asymmetry of $\Omega$ is
\[
 \mathcal A(\Omega)=\inf_{a\in\R^n}
 \frac{|\Omega\mathbin\triangle(a+\Omega^*)|}{|\Omega|},
\]
where $E\triangle F=(E\setminus F)\cup(F\setminus E)$ is the
symmetric difference. Define the normalized spectral deficits by
\[
 \delta_\Gamma(\Omega)=\frac{\Gamma_1(\Omega)}{\Gamma_1(\Omega^*)}-1,
 \quad
 \delta_\lambda(\Omega)=\frac{\lamD(\Omega)}{\lamD(\Omega^*)}-1.
\]
We suppress the argument $\Omega$ when the domain is fixed.
\begin{corollary}\label{cor:quantitative}
For $n\in\{4,5\}$ there is $\kappa_n>0$ such that
\begin{equation}\label{eq:quantitative}
 \delta_\Gamma\ge2\delta_\lambda+\delta_\lambda^2
 \ge2\kappa_n\mathcal A(\Omega)^2+
       \kappa_n^2\mathcal A(\Omega)^4.
\end{equation}
The power $2$ cannot be replaced by a smaller positive power in a uniform
lower bound for $\delta_\Gamma$ in terms of $\mathcal A$.
\end{corollary}
The lower bound follows from the sharp quantitative Faber-Krahn
inequality of Brasco, De Philippis, and Velichkov \cite{BDPV}. We give an
explicit affine-deformation argument for optimality of the power in
Section~\ref{sec:quantitative}.

The sign conditions explain the dimensional restriction in this proof.
Remark~\ref{rem:higher} shows that the same functional comparison
already fails in dimension six.

Section~\ref{sec:radial} constructs the radial transformation and
derives its basic identities. In Section~\ref{sec:functional}, we prove
the defect identity for arbitrary real trial functions.
Section~\ref{sec:global} establishes the profile inequalities needed
for positivity in dimensions four and five. The proof of
Theorem~\ref{thm:main} is completed in
Sections~\ref{sec:comparison}-\ref{sec:equality}: the first proves the
sharp spectral inequalities, and the second treats equality.
Section~\ref{sec:quantitative} proves Corollary~\ref{cor:quantitative}
and the optimality of the quadratic power of the Fraenkel asymmetry.

\section{The radial transformation}\label{sec:radial}

We first work in an arbitrary dimension $n\ge2$. Put $m=n-1$,
$\nu=n/2-1$, and suppress the dimension subscripts on $q$ and $\beta$.
A prime denotes differentiation in $r$. Use the normalized Bessel
functions
\begin{equation}\label{eq:normbessel}
 j(r)=2^\nu\Gamma(\nu+1)r^{-\nu}J_\nu(r),\quad
 i(r)=2^\nu\Gamma(\nu+1)r^{-\nu}I_\nu(r).
\end{equation}
Here $\Gamma$ is Euler's gamma function. The functions $j,i$ are
even and entire, with $j(0)=i(0)=1$. Define
\begin{gather}
 c=\frac{j(\beta)}{i(\beta)},\quad
 u=j-ci,\quad p=-u',\quad
 \phi(r)=j(qr/\beta),\quad \lambda=q^2/\beta^2,\label{eq:profiles}\\
 p_1=u''(\beta)=-2j(\beta),\quad
 \rho=\frac{I_{\nu+1}(\beta)}{I_\nu(\beta)},\quad h=\beta\rho,
 \quad K=\left(\frac{\phi'(\beta)}{p'(\beta)}\right)^2,\label{eq:K}\\
 \ell=\frac1{1-c},\quad \gamma=\lambda\ell^2-K.\label{eq:gamma}
\end{gather}
The radial Laplacian is $\Delta v=v''+mv'/r$. These functions satisfy
\begin{equation}\label{eq:radialODE}
 \Delta^2u=u,\quad -\Delta\phi=\lambda\phi,\quad
 u(\beta)=u'(\beta)=\phi(\beta)=0,
 \quad \phi'(\beta)=-\sqrt Kp_1.
\end{equation}

\begin{lemma}\label{lem:coordinate}
One has $-1<c<0$ and $p_1>0$. The functions $u,p,\phi$ are positive on
$(0,\beta)$, and $t=u/p$ is a smooth strictly decreasing bijection from
$(0,\beta)$ onto $(0,\infty)$.
\end{lemma}
\begin{proof}
Let $\tau=j_{\nu+1,1}$. The Bessel recurrence gives
\[
 j'(r)=-2^\nu\Gamma(\nu+1)r^{-\nu}J_{\nu+1}(r)<0
 \qquad(0<r<\tau).
\]
If $q\ge\tau$, integrating
$(r^{n-1}j')'=-r^{n-1}j$ would give
\[
 0=\tau^{n-1}j'(\tau)
   =-\int_0^\tau s^{n-1}j(s)\,ds<0,
\]
a contradiction. Thus $q<\tau$ and $j(\tau)<0$. The cross-product
in \eqref{eq:crossproduct} is positive on $(0,q]$ and negative at
$\tau$, so its first positive zero satisfies $q<\beta<\tau$.
Moreover,
\[
 J_{\nu+2}(r)
 =r^{-\nu-2}\int_0^r s^{\nu+2}J_{\nu+1}(s)\,ds>0
 \qquad(0<r\le\tau).
\]
Hence $q<\beta<j_{\nu+1,1}<j_{\nu+2,1}$ and $j(\beta)<0$,
while the defining series give $|j(\beta)|<i(\beta)$.

For $R=I_{\nu+1}/J_{\nu+1}$, the Bessel equations yield
\[
 [r(I_{\nu+1}'J_{\nu+1}-I_{\nu+1}J_{\nu+1}')]'
 =2rI_{\nu+1}J_{\nu+1}>0.
\]
The expression in brackets tends to zero as $r\to0$, so $R'>0$.
Since $1+cR(\beta)=0$ and $c<0$,
\[
 p=2^\nu\Gamma(\nu+1)r^{-\nu}J_{\nu+1}(1+cR)>0.
\]
It follows that $u(r)=\int_r^\beta p(s)\,ds>0$. The definition of $q$
gives $\phi>0$.

Set
\begin{equation}\label{eq:ACeta}
 A=-p'/p,\quad C=-m/r,\quad \eta=(\Delta u)'/p.
\end{equation}
The radial equations imply
\begin{equation}\label{eq:ACetaODE}
 A'=A^2+AC+\eta-C',\quad
 \eta'=(A+C)\eta+t,\quad C'=m/r^2.
\end{equation}
Moreover, $\eta=(1-cR)/(1+cR)$, so $\eta'>0$, and
\[
 (p/r)'=2^\nu\Gamma(\nu+1)r^{-\nu-1}(-J_{\nu+2}+cI_{\nu+2})<0.
\]
Thus $Ar>-1$. Near zero, $A'=r^{-2}+O(1)>0$. At a first zero of $A'$,
differentiating \eqref{eq:ACetaODE} would give
$A''=m(Ar+2)/r^3+\eta'>0$, a contradiction. Therefore $\log p$ is
strictly concave.

Now $t'=tA-1$. This is negative when $A\le0$. When $A>0$, strict
log-concavity gives
\[
 u(r)<\int_r^\infty p(r)e^{-A(r)(s-r)}\,ds=p(r)/A(r),
\]
so $t'<0$ again. Finally, $p(r)\sim(1+c)r/n$ at zero, while
$t(r)\sim(\beta-r)/2$ at the boundary. These are the required endpoint
limits.
\end{proof}

\subsection{Coefficient identities}

A subscript $t$ denotes differentiation in $t$. All coefficients are
evaluated at $r=r(t)$, the inverse of $t=u/p$.
Define
\begin{gather}
 T=t'=tA-1,\quad f=\phi/p,\quad z=f_t,\quad
 g=f-tz,\quad w=\phi'/p,\label{eq:fgzw}\\
 V=w^2-Cfw+\lambda f^2,\quad
 J=A'+(1+tC)\eta+t^2,\label{eq:VJ}\\
 G=(V-KJ)/C',\quad b=-K(A+t\eta)-CG-fw,\label{eq:fluxcoeff}\\
 a=K-g^2,\quad d=2G-tG_t,\quad
 B=\frac{K+(m-1)G}{m},\quad \D=aB-d^2/4.\label{eq:abcd}
\end{gather}
The term $-Kt\eta$ in $b$ compensates for the plate potential $-KU^2$.
The coefficient $B$ records the number of tangential directions.

For later use, we derive the identities behind the transformation.
The radial equations give
\begin{align}
 f'&=w+Af=Tz,&w'&=(A+C)w-\lambda f,&w&=-z-Ag,\label{eq:fw}\\
 V'&=(2A+C)V-C'fw,\notag\\*
 J'-(2A+C)J&=C'(A+t\eta)-\frac2mCC'.\label{eq:Jp}
\end{align}
Using $C''/C'=2C/m$ in the derivative of $G$ gives
\begin{equation}\label{eq:Gp}
 G'=b+2AG+2BC,\quad b_t=-2wz-CG_t-2K\eta.
\end{equation}
Substitution of \eqref{eq:fw}-\eqref{eq:Gp} yields
\begin{align}
 2b-tb_t-2zg&=-2aA-dC,\label{eq:kernel1}\\
 b+G_t&=-dA-2BC,\label{eq:kernel2}\\
 -Kt^2-z^2+\lambda f^2+b_t&=aA^2+dAC+BC^2.\label{eq:kernel3}
\end{align}
For example, the second identity follows by inserting
$G_t=(b+2AG+2BC)/(tA-1)$. The first and third then follow from
$w=-z-Ag$ and the definitions of $V$ and $J$.

\begin{proposition}[Profile positivity]\label{prop:signs}
For $n\in\{4,5\}$ and $0<r<\beta$,
\begin{equation}\label{eq:signs}
 a>0,\quad B>0,\quad \D>0,\quad K-G>0.
\end{equation}
Furthermore, $\gamma>3/50$ when $n=4$, and
$.0244<\gamma<.0245$ when $n=5$.
\end{proposition}
Section~\ref{sec:global} proves this proposition. It also implies the
remaining sign needed in the Hessian identity. Indeed, $G<K$ gives
$B<K$, and hence
\begin{equation}\label{eq:cross}
 2K-g^2-d=K+a-d\ge K+a-2\sqrt{aB}
 \ge(\sqrt K-\sqrt a)^2\ge0.
\end{equation}

\subsection{Endpoint behavior}

Let $\delta=\beta-r$. Taylor expansion in the radial equations gives
\begin{gather}
 t=\delta/2+O(\delta^2),\quad
 f(t)=\sqrt K(1+\rho t+O(t^2)),\label{eq:fboundary}\\
 G=\frac K3(2/\beta-\rho)\delta+O(\delta^2),\quad
 b(t)=K\rho+O(t).\label{eq:bGboundary}
\end{gather}
We use the continuous endpoint values $f(0)=\sqrt K$ and $G(0)=0$.
The value $b(0+)=K\rho>0$ determines the sign of the nodal contribution.
It can be checked directly from
\[
\begin{gathered}
 A=\delta^{-1}-(\rho-m/\beta)/2+O(\delta),\quad
t\eta=\rho/2+O(\delta),\\[3pt]
 w=-\sqrt K\{\delta^{-1}+m/\beta+(\rho-m/\beta)/2+O(\delta)\}.
\end{gathered}
\]
At the center, the even analytic expansions imply
\begin{equation}\label{eq:centerasym}
 f(t)=\ell t+O(t^{-1}),\quad
 G(t)=G_\infty+O(t^{-2}),\quad
 b(t)=b_\infty t+O(t^{-1})\quad(t\to\infty).
\end{equation}
These are convergent expansions in $t^{-1}$ near infinity. In particular,
\begin{equation}\label{eq:boundedcoeffs}
 f_t,\quad g,\quad G,\quad G_t,\quad tG_t,\quad b_t,\quad b-tb_t
\end{equation}
are bounded on $(0,\infty)$. These conclusions hold in every fixed
dimension and do not use Proposition~\ref{prop:signs}.

\section{The defect identity}\label{sec:functional}

In this section, let $n\in\{4,5\}$. Integrals without a specified
domain are over $\R^n$.
For $U\in C_c^\infty(\Omega)$, put $P=|\nabla U|$ and define
\begin{equation}\label{eq:Fextension}
 F[U]=\begin{cases}
 P f(|U|/P),&P>0,\\
 \ell|U|,&P=0.
 \end{cases}
\end{equation}
The bounds for $f_t$ and $g=f-tf_t$ show that the map
$(v,\xi)\mapsto|\xi|f(|v|/|\xi|)$, with this extension at $\xi=0$,
is globally Lipschitz. Therefore $F[U]\in H_0^1(\Omega)$ and has compact
support. It is nonzero whenever $U$ is nonzero.

\subsection{The pointwise identity}

Write $D^2U=(\partial_{ij}U)$ for the Hessian. Matrix norms are
Frobenius norms, so $|D^2U|^2=\sum_{i,j}(\partial_{ij}U)^2$.
First suppose $U>0$ and $P>0$. Choose an orthonormal frame
$(e_1,\ldots,e_n)$ with $e_1=\nabla U/P$, set
$U_{ij}=D^2U(e_i,e_j)$, and write
\[
 X=U_{11},\quad Y=\sum_{i=2}^nU_{ii},\quad
 Z=(U_{12},\ldots,U_{1n}).
\]
Let $W_0$ be the trace-free part of the tangential $m\times m$ Hessian
block. Then
\begin{equation}\label{eq:Hess}
 |D^2U|^2=X^2+Y^2/m+2|Z|^2+|W_0|^2,
 \quad |\nabla F|^2=(zP+gX)^2+g^2|Z|^2.
\end{equation}
Set $T_1(H)=(\tr H)I_n-H$, where $\tr$ is the trace and $I_n$
is the identity matrix. The identity
$\Div T_1(D^2U)=0$ and direct differentiation give
\begin{align}
 \Div\{Pb(U/P)\nabla U\}
 &=b_tP^2+(2b-tb_t)PX+bPY,\label{eq:div1}\\
 \Div\{G(U/P)T_1(D^2U)\nabla U\}
 &=dXY+\frac{m-1}{m}GY^2-d|Z|^2-G|W_0|^2+G_tPY.
 \label{eq:div2}
\end{align}
Combining these formulas with \eqref{eq:kernel1}-\eqref{eq:kernel3}
proves
\begin{align}
 &K(|D^2U|^2-U^2)-(|\nabla F|^2-\lambda F^2)
   +\Div\{Pb\nabla U+GT_1(D^2U)\nabla U\}\notag\\*
 &\quad=a(X-AP)^2+d(X-AP)(Y-CP)+B(Y-CP)^2\notag\\*
 &\quad\quad +(2K-g^2-d)|Z|^2+(K-G)|W_0|^2.
 \label{eq:identity}
\end{align}
All coefficients are evaluated at $t=U/P$. By
Proposition~\ref{prop:signs} and \eqref{eq:cross}, the right side is
nonnegative. For $U<0$, apply the same formula to $-U$.

\subsection{Integration across nodal and critical sets}\label{sec:weak}

Extend $f,G$ evenly and $b$ oddly away from zero. With $b_0=K\rho>0$,
put
\[
 b_{\mathrm{reg}}(t)=b(t)-b_0\operatorname{sign}t,
 \quad b_{\mathrm{reg}}(0)=0.
\]
The endpoint expansions and \eqref{eq:boundedcoeffs} imply that
$P b_{\mathrm{reg}}(U/P)$ is a Lipschitz function of $(U,\nabla U)$,
with extension $b_\infty U$ at $P=0$. Likewise,
$G(U/P)\nabla U$ is Lipschitz, with value zero at $P=0$.
Consequently the vector fields
\[
 V_1=P b_{\mathrm{reg}}(U/P)\nabla U,\quad
 V_2=G(U/P)T_1(D^2U)\nabla U
\]
are compactly supported and Lipschitz. Their divergences have integral
zero.

For $W=P\nabla U$, which is $C^1$ and compactly supported, and
$\theta_\varepsilon(s)=s/(s^2+\varepsilon^2)^{1/2}$, integration by
parts gives
\[
 \int\theta_\varepsilon(U)\Div W
 =-\int\theta_\varepsilon'(U)P^3\le0.
\]
Dominated convergence on the left proves
\begin{equation}\label{eq:nodalsign}
 \int\operatorname{sign}(U)\Div W\le0.
\end{equation}
No assumption on the number or area of the nodal components is involved.

On $\{U\ne0,P>0\}$, the pointwise identity now reads
\begin{equation}\label{eq:weakdensity}
 D_U+\Div V_1+\Div V_2+b_0\operatorname{sign}(U)\Div W=\Q_U,
\end{equation}
where
$D_U=K(|D^2U|^2-U^2)-(|\nabla F|^2-\lambda F^2)$ and $\Q_U$ is the
right side of \eqref{eq:identity}, applied to $|U|$ with normal direction
$\nabla|U|/P$. The set $\{U=0,P>0\}$ is locally a smooth hypersurface
and is Lebesgue null. On $\{P=0\}$, locality of weak derivatives gives
$D^2U=0$ almost everywhere. Since $F=\ell|U|$ there,
$\nabla F=0$ almost everywhere on that set as well. The vector fields
$V_1,V_2,W$ and their weak derivatives vanish there almost everywhere.
Thus the remaining density is $D_U=\gamma U^2\ge0$.
Integrating \eqref{eq:weakdensity} proves
\begin{equation}\label{eq:functional}
 \int|\nabla F|^2-\lambda\int F^2
 \le K\left(\int(\Delta U)^2-\int U^2\right),
\end{equation}
where $\int|D^2U|^2=\int(\Delta U)^2$ by integration by parts.

The same argument gives the full remainder.
\begin{proposition}\label{prop:defect}
For every real $U\in C_c^\infty(\Omega)$, set
$N_U=\{U=0,\,|\nabla U|>0\}$. Then
\begin{align}
 &K\left(\int(\Delta U)^2-\int U^2\right)
 -\left(\int|\nabla F[U]|^2-\lambda\int F[U]^2\right)\notag\\*
 &\quad=\int_{\{U\ne0,P>0\}}\Q_U
 +\gamma\int_{\{P=0\}}U^2
 +2K\rho\int_{N_U}P^2\,d\mathcal H^{n-1}.\label{eq:defect}
\end{align}
Every term on the right is finite and nonnegative.
\end{proposition}
\begin{proof}
It remains to identify the nodal term in \eqref{eq:nodalsign}.
Choose a smooth increasing function $\chi:[0,\infty)\to[0,1]$ equal
to zero on $[0,1/2]$ and one on $[1,\infty)$. The field
$W\chi(P/\delta)$ vanishes near the critical set. On its support the
nodal set is a smooth hypersurface, so the distributional jump formula
and integration give
\begin{align*}
 2\int_{N_U}\chi(P/\delta)P^2\,d\mathcal H^{n-1}
 ={}&-\int\operatorname{sign}(U)\chi(P/\delta)\Div W\\
 &-\int\operatorname{sign}(U)W\cdot\nabla\chi(P/\delta).
\end{align*}
The last integral is bounded in absolute value by
$C_*\delta\int_{\supp\nabla U}|D^2U|$, where $C_*$ is
independent of $\delta$. It therefore tends to zero.
Also $\Div W=0$ almost everywhere on $\{P=0\}$. Dominated convergence
on the right and monotone convergence on the left yield
\[
 -\int\operatorname{sign}(U)\Div W
 =2\int_{N_U}P^2\,d\mathcal H^{n-1}<\infty.
\]
Integrating \eqref{eq:weakdensity}, including the density $\gamma U^2$
on the critical set, proves \eqref{eq:defect}.
\end{proof}

\section{Global positivity in dimensions four and five}\label{sec:global}

This section proves Proposition~\ref{prop:signs}. Set
\[
 x=1-r^2/\beta^2,\quad y=1-x,\quad
 Y(x)=\frac{8u(r)}{p_1\beta^2},\quad
 Z(x)=\frac{2\phi(r)}{\sqrt Kp_1\beta}.
\]
The boundary corresponds to $x=0$ and the center to $x=1$.
The even Bessel series show that $Y,Z$ are entire in $x$.
Within this section, primes on the entire profiles in $x$ denote
$x$ derivatives. Derivatives of the original radial functions remain
in $r$.
For an analytic function $A$, write $A_k=A^{(k)}(0)/k!$ for the
coefficient of $x^k$.
Their equations and initial coefficients are
\begin{align}
 y^2Y''''-(n+2)yY'''+\frac{n(n+2)}4Y''&=\frac{\beta^4}{16}Y,
 \label{eq:YZode}\\
 yZ''-\frac n2Z'+\frac{q^2}4Z&=0,\notag\\*
 Y_0=Y_1=0,\quad Y_2=1,\quad Y_3&=(n+2-h)/6,
 \quad Z_0=0,\quad Z_1=1,\quad Z_2=n/4.\label{eq:YZinitial}
\end{align}

\subsection{Entire expressions for the coefficients}

Define
\begin{gather}
 L=Y/x^2,\quad P=Y'/(2x),\quad Q=P+xP',\quad
 \Phi=Z/x,\quad\Psi=Z',\label{eq:LP}\\
 W=4yP^2+xLP-2yLQ,\quad H=2P\Phi-L\Psi.\label{eq:WH}
\end{gather}
The symbols $Y,Z,P$ in this section denote entire functions of $x$.
Put
\begin{align}
 \mathscr J={}&(Y'-2yY'')^2+2yY'(3Y''-2yY''')\notag\\*
 &+(4yY'-2mY)\left(\frac{n+2}{2}Y''-yY'''\right)
       +\frac{\beta^4}{4}Y^2,\label{eq:Jscript}\\
 N={}&16y(Z')^2-8mZZ'+4q^2Z^2-\mathscr J,\quad
 \calR=-N/x^3,\label{eq:Rscript}\\
 \calA={}&(W^2-4yH^2)/x^2,\quad
 \calB=4mP^2-(m-1)x\calR,\label{eq:ABscript}\\
 \calD={}&\calR W-yL\big[(\calR+x\calR')P-2x\calR P'\big],
 \label{eq:Dscript}\\
 \calE={}&4P^2\calA\calB-\calD^2.\label{eq:Escript}
\end{align}
Equations~\eqref{eq:YZode}-\eqref{eq:YZinitial} give
$N_0=N_1=N_2=0$ and
$(W^2-4yH^2)_0=(W^2-4yH^2)_1=0$.
All apparent divisions above are therefore removable, and every displayed
function is entire.

Direct differentiation gives
\[
 t=\frac{\beta xL}{4\sqrt yP},\quad
 \frac{dt}{dr}=-\frac{W}{4yP^2}.
\]
Thus $P,W>0$ on $(0,1)$ by Lemma~\ref{lem:coordinate}.
At the endpoints, $W(0)=2$ and $W(1)=L(1)P(1)>0$.
Substitution in \eqref{eq:abcd} gives
\begin{gather}
 \frac g{\sqrt K}=\frac{2\sqrt yH}{W},\quad
 \frac GK=-\frac{x\calR}{4mP^2},\label{eq:xgG}\\
 \frac aK=\frac{x^2\calA}{W^2},\quad
 \frac BK=\frac{\calB}{4m^2P^2},\quad
 \frac dK=-\frac{x\calD}{2mP^2W},\label{eq:xabd}\\
 \frac{\D}{K^2}=\frac{x^2\calE}{16m^2P^4W^2}.\label{eq:xdet}
\end{gather}
It remains to prove the following bounds on $0\le x\le1$:
\begin{equation}\label{eq:envelopes}
\begin{array}{c|cc}
 &n=4&n=5\\ \hline
 P&\displaystyle\ge\frac56(1+x)&\displaystyle\ge\frac{99}{100}(1+x)\\[3pt]
 \calR&\ge1&>3\\[2pt]
 \calA&\displaystyle\ge\frac35(1+x)^4&\displaystyle\ge\frac25(1+x)^6\\[3pt]
 \calB&\displaystyle\ge6(1+x)^2&\displaystyle\ge\frac{54}{5}(1+x)^2\\[3pt]
 |\calD|&\le3(1+x)^4&\displaystyle\le\frac{81}{20}(1+x)^5.
\end{array}
\end{equation}
They imply
\begin{equation}\label{eq:Epositive}
 \calE\ge
 \begin{cases}
 (1+x)^8,&n=4,\\[2pt]
 \dfrac{133407}{250000}(1+x)^{10},&n=5,
 \end{cases}
\end{equation}
because
\[
 4(5/6)^2(3/5)6-3^2=1,\quad
 4(99/100)^2(2/5)(54/5)-(81/20)^2=133407/250000.
\]
These bounds leave positive margins in both dimensions.

\subsection{Scalar bounds and Taylor coefficients}\label{sec:scalar}

Put $b_n=-j'/r$ and $b_{n,I}=i'/r$, where the derivatives here are
in $r$. Write $(a)_0=1$ and $(a)_k=\prod_{j=0}^{k-1}(a+j)$ for
the rising factorial. The defining series are
\begin{equation}\label{eq:scalarseries}
 j(r),i(r)=\sum_{k\ge0}\frac{(\mp1)^k(r^2/4)^k}{k!(n/2)_k},
 \quad
 b_n(r),b_{n,I}(r)=\frac1n\sum_{k\ge0}
 \frac{(\mp1)^k(r^2/4)^k}{k!(n/2+1)_k}.
\end{equation}
For $n=4$ and $0\le r\le5$, the remainders after index $19$ are
bounded in absolute value by
\[
 \frac{2(25/4)^{20}}{20!21!},\quad
 \frac{(25/4)^{20}}{20!22!},
\]
respectively. For $n=5$ and $0\le r\le6$, the corresponding bounds
after index $23$ are
\[
 \frac{2\,9^{24}}{24!(5/2)_{24}},\quad
 \frac{2\,9^{24}}{5\cdot24!(7/2)_{24}}.
\]
In each case the ratios of successive omitted absolute terms are less
than $1/2$.

The finite sums with these remainder bounds give the following strict
enclosures. Every decimal in the proof denotes an exact rational number.
\begingroup\small
\begin{equation}\label{eq:rootbounds}
\begin{array}{c|cc}
 &n=4&n=5\\\hline
 q&(3.831705970207512,\ 3.831705970207513)
   &(4.493409457909064,\ 4.493409457909065)\\
 \beta&(4.610899879049055,\ 4.610899879049056)
   &(5.267657530336814,\ 5.267657530336815)
\end{array}
\end{equation}
\endgroup
At the endpoints the tested expressions are $j$ and $b_ni+jb_{n,I}$,
respectively; both are positive on the left and negative on the right.
For example, in dimension four the two $j$ values are greater than
$6\cdot10^{-17}$ and less than $-14\cdot10^{-17}$, and the two
cross-product values are greater than $29\cdot10^{-17}$ and less than
$-6\cdot10^{-17}$.

The following elementary bounds identify these as first roots, without
requiring numerical information about other zeros. In dimension four,
put $v=1-r^2/25$. The alternating series through index $5$ gives
\[
 4b_4(r)\ge
 \frac{138379+3887925v+4178750v^2+3556250v^3
       +234375v^4+390625v^5}{12386304}>0
\]
for $0\le r\le5$. In dimension five, with $v=1-r^2/28$, the same
truncation gives
\[
 5b_5(r)\ge
 \frac{9604v^5+3430v^4+81340v^3+89810v^2+93770v+11621}{289575}>0
 \quad(0\le r\le\sqrt{28}).
\]
Hence $j'=-rb_n<0$ up to the respective endpoint, identifying $q$.
For $D_\nu=J_\nu I_{\nu+1}+J_{\nu+1}I_\nu$,
\[
 (rD_\nu)'=2rJ_\nu I_\nu.
\]
Thus $rD_\nu$ increases before $q$ and decreases between $q$ and that
endpoint. This identifies the bracketed root as $\beta$.

Substitution in the same series yields
\begingroup\small
\begin{equation}\label{eq:parameterbounds}
\begin{array}{c|cc}
 &n=4&n=5\\\hline
 q^2&(14.681970642123,\ 14.681970642125)
     &(20.190728556426,\ 20.190728556428)\\
 \beta^4&(452.004510133173,\ 452.004510133175)
     &(769.963483241901,\ 769.963483241903)\\
 h&(3.218890982799,\ 3.218890982801)
     &(3.501551477382,\ 3.501551477384).
\end{array}
\end{equation}
\endgroup
We also use the coarser consequences
\begin{align}
 n=4:\quad&21<\beta^2<22,\quad p_1>1/5,\quad -2/125<c<0,
 \quad J_2(q)>2/5,\notag\\*
 &3/5<K<121/200,\quad69/100<\lambda<7/10;\label{eq:coarse4}\\
 n=5:\quad&27.7<\beta^2<28,\quad p_1>.1488,\quad -.0088<c<0,
 \quad q^2<20.2,\notag\\*
 &q^2b_5(q)>.65,\quad .0244<\gamma<.0245.\label{eq:coarse5}
\end{align}
For the last bounds on $K$ and $\gamma$, the scalar expressions are
\[
 K=\left(\frac{q^2b_n(q)}{2\beta j(\beta)}\right)^2,
 \quad \gamma=\frac{q^2}{\beta^2(1-c)^2}-K.
\]

Equation~\eqref{eq:YZode} gives the recurrences
\begin{align}
 Y_{k+4}={}&\frac{2k+n+2}{k+4}Y_{k+3}
 -\frac{(k+n/2)(k+(n+2)/2)}{(k+3)(k+4)}Y_{k+2}\notag\\*
 &+\frac{\beta^4Y_k}{16(k+1)(k+2)(k+3)(k+4)},\label{eq:Yrec}\\
 Z_{k+2}={}&\frac{k+n/2}{k+2}Z_{k+1}
 -\frac{q^2Z_k}{4(k+1)(k+2)}.\label{eq:Zrec}
\end{align}
Together with \eqref{eq:YZinitial} and \eqref{eq:parameterbounds}, they
give Tables~\ref{tab:coeff4} and \ref{tab:coeff5}. The first three
columns are lower bounds and the last is an upper bound on the absolute
value. A dash indicates an unused coefficient; the constant term of
$\calB$ is exact.

\begin{table}[htbp]
\centering\small
\begin{tabular}{rcccc}
\toprule
$k$&$\calR_k$&$\calA_k$&$\calB_k$&$|\calD_k|$\\
\midrule
0&2&.77&12&2.5\\
1&4&2.88&11.8&5.1\\
2&2&3.53&7.08&8.9\\
3&.9&2.90&3.78&9.3\\
4&.29&1.94&1.66&6.77\\
5&.069&1.05&.52&3.64\\
6&.0126&.47&.138&1.5\\
7&.0019&.18&.028&.5\\
8&.00023&.058&.0048&.14\\
9&---&.016&---&.04\\
10&---&.0039&---&.01\\
11&---&---&---&.01\\
12&---&---&---&.02\\
13&---&---&---&.06\\
\bottomrule
\end{tabular}
\caption{Rational coefficient bounds in dimension four.}\label{tab:coeff4}
\end{table}

\begin{table}[htbp]
\centering\small
\begin{tabular}{rcccc}
\toprule
$k$&$\calR_k$&$\calA_k$&$\calB_k$&$|\calD_k|$\\
\midrule
0&4&.69&16&4.01\\
1&8.14&3.53&15.97&10.51\\
2&6.16&5.65&6.45&21.10\\
3&3.216&5.87&3.01&27.16\\
4&1.199&4.84&1.53&24.990\\
5&.3516&3.26&.65&17.248\\
6&.0810&1.82&.23&9.327\\
7&.0155&.85&.0702&4.113\\
8&.0024&.347&.0168&1.513\\
9&.00026&.121&.0032&.473\\
10&---&.0375&.00030&.129\\
11&---&.0102&$-.000604$&.043\\
12&---&.00247&$-.001940$&.050\\
13&---&.00050&$-.005538$&.137\\
14&---&$-.000006$&---&.405\\
15&---&$-.000274$&---&1.192\\
\bottomrule
\end{tabular}
\caption{Rational coefficient bounds in dimension five.}\label{tab:coeff5}
\end{table}

The finite arithmetic in these tables uses only
\begin{equation}\label{eq:coefficientoperations}
 (AB)_k=\sum_{j=0}^kA_jB_{k-j},\quad
 (A')_k=(k+1)A_{k+1},\quad (xA)_k=A_{k-1},
\end{equation}
with negative indices zero, applied to \eqref{eq:LP}-\eqref{eq:Dscript}
in that order. Every parameter is enclosed in its rational interval in
\eqref{eq:parameterbounds}; the divisions by $x$ shift coefficients.
Thus the tables involve no estimates of a quotient near a removable
singularity. The upper bounds for $|\calD_k|$, $4\le k\le13$, in
dimension four sum to $12.69<13$. In dimension five they imply
\begin{equation}\label{eq:finiteMass}
 \sum_{k=8}^{15}\max\{-\calA_k,0\}<.0003,\quad
 \sum_{k=7}^{13}\max\{-\calB_k,0\}<.0081,\quad
 \sum_{k=4}^{15}|\calD_k|<60.
\end{equation}

\subsection{Bounds for the full infinite remainders}

We record explicit majorants to control every omitted coefficient.
Differentiating the Bessel series in $r^2$ gives, for $k\ge2$,
\begin{align}
 Y_k&=\frac{8\,2^\nu\Gamma(\nu+1)\beta^{k-\nu-2}}
 {p_1\,2^kk!}
 \big[J_{\nu+k}(\beta)-c(-1)^kI_{\nu+k}(\beta)\big],
 \label{eq:Ybesselcoef}\\
 Z_k&=\frac{q^{k-1}J_{\nu+k}(q)}{2^{k-1}k!J_{\nu+1}(q)}.
 \label{eq:Zbesselcoef}
\end{align}
For $k\ge4$ in dimension four and $k\ge5$ in dimension five, the
normalized $J$ series alternates with decreasing term magnitudes.
Its absolute value is at most one. In these ranges, the normalized $I$
series is bounded by \mbox{$\exp(\beta^2/[4(\nu+k+1)])<3$}. Thus
\begin{align}
 n=4:\quad |Y_k|&\le\frac{2(11/2)^k}{k!(k+1)!},&
 |Z_k|&\le\frac{(5/2)(15/4)^k}{k!(k+1)!}\quad(k\ge4),
 \label{eq:factorial4}\\
 n=5:\quad |Y_k|&\le\frac{2\cdot7^k}{k!(5/2)_k},&
 |Z_k|&\le\frac{(31/10)(101/20)^k}{k!(5/2)_k}\quad(k\ge5).
 \label{eq:factorial5}
\end{align}
For example, in dimension five the prefactors follow from
$8(1+3|c|)/(p_1\beta^2)<2$ and $2/(q^2b_5(q))<31/10$.
At $k=11$, the bounds in dimension four are less than
$3/(2\cdot10^8)$ and $27/10^{11}$; in dimension five they are less than
$8/10^8$ and $33/10^{10}$. Every subsequent ratio is less than $1/20$.

Write $A\preceq A^+$ when $|A_k|\le A_k^+$ for all $k$, with
$A_k^+\ge0$. The recurrences give the finite bounds in the following
majorants, and \eqref{eq:factorial4}-\eqref{eq:factorial5} give their
tails:
\begin{align}
 Y_4^+(x)={}&x^2+.465x^3+.196x^4+.035x^5+.0056x^6+.00052x^7
 \notag\\*
 &+.000045x^8+.0000026x^9+.00000014x^{10}
 +\frac{3x^{11}}{2\cdot10^8(1-x/20)},\label{eq:Ymajor4}\\
 Z_4^+(x)={}&x+x^2+.389x^3+.083x^4+.0112x^5+.00105x^6
 +.000073x^7\notag\\*
 &+.0000039x^8+.000000164x^9+.0000000056x^{10}
 +\frac{27x^{11}}{10^{11}(1-x/20)},\label{eq:Zmajor4}\\
 Y_5^+(x)={}&x^2+.584x^3+.292x^4+.0651x^5+.0127x^6+.0015x^7
 \notag\\*
 &+.000159x^8+.0000116x^9+.00000078x^{10}
 +\frac{8x^{11}}{10^8(1-x/20)},\label{eq:Ymajor5}\\
 Z_5^+(x)={}&x+\tfrac54x^2+.618x^3+.169x^4+.0295x^5+.00362x^6
 +.000329x^7
 \notag\\*
 &+.0000231x^8+.00000129x^9+.0000000586x^{10}
 +\frac{33x^{11}}{10^{10}(1-x/20)}.\label{eq:Zmajor5}
\end{align}
The subscripts $4,5$ on these majorants indicate dimension. Thus
$Y\preceq Y_n^+$ and $Z\preceq Z_n^+$.

To propagate the majorants, write $Y^+=Y_n^+$, $Z^+=Z_n^+$,
$v=1+x$, and define the coefficient shift
\[
 \mathcal T_sA=\frac{A-\sum_{k<s}A_kx^k}{x^s}.
\]
Let $(Q_n,H_n)=(15,453)$ for $n=4$ and $(101/5,770)$ for $n=5$;
these bound $q^2,\beta^4$ from above. Put
\begin{align*}
 \mathscr J^+={}&((Y^+)'+2v(Y^+)'')^2
 +2v(Y^+)'(3(Y^+)''+2v(Y^+)''')\\
 &+(4v(Y^+)'+2mY^+)
       \left(\frac{n+2}{2}(Y^+)''+v(Y^+)'''\right)
       +\frac{H_n}{4}(Y^+)^2,\\
 N^+={}&16v((Z^+)')^2+8mZ^+(Z^+)'+4Q_n(Z^+)^2+\mathscr J^+,\\
 \calR^+={}&\mathcal T_3N^+,\quad
 L^+=\mathcal T_2Y^+,\quad P^+=\tfrac12\mathcal T_1(Y^+)',\\
 Q^+={}&P^++x(P^+)',\quad
 \Phi^+=\mathcal T_1Z^+,\quad\Psi^+=(Z^+)',\\
 W^+={}&4v(P^+)^2+xL^+P^++2vL^+Q^+,\quad
 H^+=2P^+\Phi^++L^+\Psi^+,\\
 \calA^+={}&\mathcal T_2\big[(W^+)^2+4v(H^+)^2\big],\quad
 \calB^+=4m(P^+)^2+(m-1)x\calR^+,\\
 \calD^+={}&\calR^+W^+
 +vL^+\big[(\calR^++x(\calR^+)')P^++2x\calR^+(P^+)'\big].
\end{align*}
Addition, multiplication, differentiation, and coefficient shifting
preserve majorization. Hence
$\calR\preceq\calR^+$, $\calA\preceq\calA^+$,
$\calB\preceq\calB^+$, and $\calD\preceq\calD^+$.
Every majorant is an explicit rational function with nonnegative
coefficients and converges at one. For any such pair,
\begin{equation}\label{eq:tailrule}
 \sum_{k\ge s}|A_k|\le A^+(1)-\sum_{k<s}A_k^+.
\end{equation}
Evaluating the rational expressions gives
\begin{equation}\label{eq:tailtable}
\begin{array}{c|cc|cc}
 &\multicolumn{2}{c|}{n=4}&\multicolumn{2}{c}{n=5}\\
 A&s&\text{tail bound}&s&\text{tail bound}\\\hline
 \calR&9&1/20&10&.064\\
 \calA&11&5/4&16&.008\\
 \calB&9&4/5&14&.0003\\
 \calD&14&9/2&16&4.5
\end{array}
\end{equation}
All bounds are strict.

\subsection{Proof of the envelopes}

First, $Y_k>0$ for all $k\ge2$. Indeed, $J_{\nu+1}>0$ on
$(0,\beta]$, and
\[
 J_{a+1}(z)=z\int_0^1s^{a+1}J_a(zs)\,ds
\]
propagates this positivity to all the required larger orders.
The same identity holds for $I_a$. The Wronskian relation shows that
$J_a/I_a$ decreases strictly in $z$ while $J_a>0$; the integral
identities then give
\[
 \frac{J_{a+1}(z)}{I_{a+1}(z)}>
 \frac{J_a(z)}{I_a(z)},\quad 0<z\le\beta.
\]
Since $J_{\nu+1}(\beta)/I_{\nu+1}(\beta)=-c$,
\eqref{eq:Ybesselcoef} proves the assertion. Thus every coefficient of
$P$ is positive. In dimension four,
$P_0=1$ and $P_1=(6-h)/4>2/3$, whence
\[
 P\ge1+\frac23x\ge\frac56(1+x).
\]
In dimension five, the recurrences give $P_1>.87$ and $P_2>.58$.
Therefore $P\ge1+.87x+.58x^2\ge.99(1+x)$, because
\[
 .01-.12x+.58x^2
 =\frac{29}{50}\left(x-\frac3{29}\right)^2+\frac{11}{2900}>0.
\]

\paragraph{Dimension four.}
Table~\ref{tab:coeff4} and \eqref{eq:tailtable} imply
\begin{align*}
 \calR&>2-x^9\ge1,\\
 \calA&>.77+2.88x+3.53x^2+2.90x^3+1.94x^4-2x^{11},\\
 \calB&\ge12+10x+6x^2-2x^9,\\
 |\calD|&<2.5+5.1x+8.9x^2+9.3x^3+13x^4+5x^{14}.
\end{align*}
After subtracting $\frac35(1+x)^4$ from the lower bound for $\calA$
and using $x^{11}\le x^4$, the difference is at least
\[
 .17+.48x-.07x^2+.50x^3-.66x^4\ge .17+.25x>0.
\]
Also $\calB\ge12+8x+6x^2\ge6(1+x)^2$.
For $\calD$, use $13x^4+5x^{14}\le18x^4$. The difference from
$3(1+x)^4$ is then at least
\[
 .5+6.9x+9.1x^2+2.7x^3-15x^4\ge.5+x+2.7x^3>0.
\]
This proves the dimension-four column of \eqref{eq:envelopes}.

\paragraph{Dimension five.}
Here $\calR>4-.064x^{10}>3$.
Table~\ref{tab:coeff5}, \eqref{eq:finiteMass}, and \eqref{eq:tailtable}
give
\begin{align*}
 \calA\ge A_-:={}&.69+3.53x+5.65x^2+5.87x^3+4.84x^4
       +3.26x^5+1.82x^6+.84x^7,\\
 \calB\ge B_-:={}&16+15.97x+6.45x^2+3.01x^3+1.53x^4
       +.65x^5+.22x^6,\\
 |\calD|\le D_+:={}&4.01+10.51x+21.10x^2+27.16x^3+65x^4.
\end{align*}
For the first inequality, the sum of possible negative contributions
from index eight onwards is less than $.0003+.008<.01$; absorb them
by lowering the seventh coefficient from $.85$ to $.84$.
For the second, use $.0081+.0003<.01$ and lower the sixth coefficient
from $.23$ to $.22$. For the third, the remaining absolute sum is less
than $60+4.5<65$.

We check the three polynomial comparisons in the Bernstein basis.
For $p(x)=\sum_{j=0}^d a_jx^j$ and $N\ge d$,
\begin{equation}\label{eq:Bernstein}
 p(x)=\sum_{k=0}^N b_k\binom Nkx^k(1-x)^{N-k},\quad
 b_k=\sum_{j=0}^{\min(k,d)}a_j\frac{\binom kj}{\binom Nj}.
\end{equation}
For $A_- - \frac25(1+x)^6$, at degree $N=10$, the successive $b_k$
are bounded below by
\[
 \frac1{1000}(290,403,508,587,618,575,440,229,32,90,900).
\]
For $B_- - \frac{54}{5}(1+x)^2$, at degree $N=9$, they are bounded
below by
\[
 \frac1{1000}(5200,4574,3828,2996,2128,1288,566,89,31,630).
\]
For $\frac{81}{20}(1+x)^5-D_+$, at degree $N=5$, they are
\[
 \left(\frac1{25},\frac{497}{250},\frac{1469}{250},
       \frac{6519}{500},\frac{7929}{500},\frac{91}{50}\right).
\]
All entries are positive, proving the inequalities on the whole interval.
This completes \eqref{eq:envelopes} and hence \eqref{eq:Epositive}.

Finally, \eqref{eq:xgG}-\eqref{eq:xdet} imply
$a,B,\D,K-G>0$. In dimension four, \eqref{eq:coarse4} gives
\[
 \gamma>\frac{69/100}{(127/125)^2}-\frac{121}{200}>\frac3{50}.
\]
The dimension-five bound for $\gamma$ is \eqref{eq:coarse5}.
Proposition~\ref{prop:signs} is proved.

\section{Spectral consequences}

\subsection{The sharp comparison and Rayleigh's inequality}\label{sec:comparison}

Let $0\ne U\in C_c^\infty(\Omega)$ and put
\[
 Q_U=\frac{\int(\Delta U)^2}{\int U^2},\quad
 s=Q_U^{1/4},\quad\widetilde\Omega=s\Omega,
 \quad\widetilde U(x)=U(x/s).
\]
The plate quotient of $\widetilde U$ is one. Applying
\eqref{eq:functional} to this function gives a nonzero
$\widetilde F\in H_0^1(\widetilde\Omega)$ such that
\[
 \int|\nabla\widetilde F|^2\le\lambda\int\widetilde F^2.
\]
It follows that
\[
 s^{-2}\lamD(\Omega)=\lamD(\widetilde\Omega)\le\lambda,
 \quad Q_U\ge\lambda^{-2}\lamD(\Omega)^2.
\]
Taking the infimum over $C_c^\infty(\Omega)$ proves
\eqref{eq:comparison}. The Faber-Krahn inequality on finite-volume
open sets \cite[Section~1.1]{BDPV} gives
\[
 \lamD(\Omega)\ge q_n^2
 \left(\frac{\omega_n}{|\Omega|}\right)^{2/n}.
\]
Since $\lambda=q_n^2/\beta_n^2$, this proves \eqref{eq:Rayleigh}.
On $B_{\beta_n}$ the explicit clamped profile $u$ has quotient one,
while \eqref{eq:comparison} gives the reverse bound. Thus
$\Gamma_1(B_R)=\beta_n^4/R^4$, establishing sharpness of both constants.

\subsection{Equality and removable sets}\label{sec:equality}

We give the equality argument for all finite-volume open sets. First,
the infimum defining $\Gamma_1$ is attained. An $L^2$-normalized
minimizing sequence is bounded in $H^2(\R^n)$: integration by parts
controls its gradient, and $\|D^2U\|_2=\|\Delta U\|_2$.
The Sobolev embedding $H^1(\R^n)\hookrightarrow L^{2n/(n-2)}(\R^n)$
gives the uniform tail estimate
\[
 \|U\|_{L^2(\Omega\setminus B_R)}^2
 \le|\Omega\setminus B_R|^{2/n}\|U\|_{2n/(n-2)}^2\longrightarrow0.
\]
Local Rellich compactness therefore gives strong $L^2$ convergence of
a subsequence. Weak closedness of $H_0^2(\Omega)$ and lower
semicontinuity give a minimizer.

We also need simplicity of the first ball eigenvalue. Here is a direct
argument in the dimensions under consideration. For a spherical harmonic
component of degree $l$, write its radial factor as $R_l$ and put
$w=r^{(n-1)/2}R_l$. Its mass and plate energy are
\[
 \int_0^\beta w^2\,dr,\quad
 E_{\mu_l}[w]=\int_0^\beta
 \left|-w''+\frac{\mu_l}{r^2}w\right|^2\,dr,
 \quad \mu_l=\left(l+\frac{n-2}{2}\right)^2-\frac14.
\]
For $l\ge1$, integration by parts yields
\begin{align*}
 E_{\mu_l}[w]-E_{\mu_0}[w]
 = (\mu_l-\mu_0)\left\{
 2\int_0^\beta\frac{w'^2}{r^2}\,dr
 +(\mu_l+\mu_0-6)\int_0^\beta\frac{w^2}{r^4}\,dr\right\}.
\end{align*}
Expanding the square in
\mbox{$\int_0^\beta(w'-3w/(2r))^2/r^2\,dr\ge0$} gives the
weighted Hardy inequality
\mbox{$\int_0^\beta w'^2/r^2\,dr\ge(9/4)\int_0^\beta w^2/r^4\,dr$}.
Consequently,
\begin{equation}\label{eq:angular}
 E_{\mu_l}[w]-E_{\mu_0}[w]
 \ge(\mu_l-\mu_0)(\mu_l+\mu_0-3/2)
       \int_0^\beta\frac{w^2}{r^4}\,dr>0
\end{equation}
for a nonzero component, since $n\ge4$.
For a smooth ball eigenfunction, $R_l=O(r^l)$ at zero. Its radial
replacement $R_l(|x|)$ belongs to $H_0^2(B_\beta)$; cutting it off
near zero has squared $H^2$ cost $O(\varepsilon^{2l+n-4})$.
Thus $E_{\mu_0}[w]\ge\int_0^\beta w^2\,dr$ by
$\Gamma_1(B_\beta)=1$, and
\eqref{eq:angular} excludes every $l\ge1$ component from a first
eigenfunction. The regular radial solutions of $\Delta^2v=v$ are
spanned by $j,i$; their two clamped conditions have rank one at $\beta$.
Every first ball eigenfunction is therefore a multiple of the positive
profile $u$.

Suppose equality holds in \eqref{eq:Rayleigh}. Equality must then hold
in the Faber-Krahn inequality, so $\Omega$ agrees almost everywhere with a translated
ball $B$. Openness gives $\Omega\subset B$: any point outside $B$,
including a point on its boundary, would have a neighborhood in $\Omega$
of positive measure outside $B$. Put $E=B\setminus\Omega$.
A minimizing eigenfunction on $\Omega$, extended by zero, is a first
eigenfunction on $B$ and hence a nonzero multiple of its positive
radial mode $u_B$.

For a compact set $K\subset\R^n$ and $k=1,2$, define the Sobolev
capacity
\[
 \operatorname{Cap}_{k,2}(K)
 =\inf\left\{\|\eta\|_{H^k(\R^n)}^2:
 \eta\in C_c^\infty(\R^n),\ \eta\ge1\text{ near }K\right\}.
\]
For open sets, take the supremum over compact subsets. For arbitrary
sets, take the infimum over open supersets. Take
$U_j\in C_c^\infty(\Omega)$ converging to $u_B$ in $H^2$.
For compact $K\subset E$, choose $\chi\in C_c^\infty(B)$ equal
to one near $K$. Then
\[
 \eta_j=\frac\chi{u_B}(u_B-U_j)
\]
is smooth, compactly supported, and equal to one near $K$, with
$\|\eta_j\|_{H^2}\to0$. Thus $\capc(K)=0$. Exhausting the relatively
closed set $E$ by compact subsets gives $\capc(E)=0$.

Conversely, let $E$ be relatively closed in $B$ with $\capc(E)=0$.
For $w\in H^2(\R^n)$, membership in $H_0^2(\Omega)$ is equivalent
to two vanishing conditions outside $\Omega$. The
$H^2$-quasicontinuous representative of $w$ must vanish there except
on a set of zero $\operatorname{Cap}_{2,2}$ capacity. The
$H^1$-quasicontinuous representatives of its first derivatives must
vanish there except on a set of zero $\operatorname{Cap}_{1,2}$
capacity. Here an $H^k$-quasicontinuous representative is continuous
outside open sets of arbitrarily small $\operatorname{Cap}_{k,2}$
capacity; see \cite{GLM24} and
\cite[Section~2, equation~(8)]{L25existence}.
Since zero $\operatorname{Cap}_{2,2}$ capacity implies zero
$\operatorname{Cap}_{1,2}$ capacity, that characterization gives
\[
 H_0^2(B\setminus E)=H_0^2(B).
\]
The eigenvalues are therefore equal. This argument applies in both
dimensions and proves the stated classification.

For a bounded $C^1$ domain, a boundary point strictly inside $B$ would
have a complement of positive measure in a small neighborhood, by the
local boundary graph. Thus equality in that class forces $\Omega=B$.
The capacity qualification is necessary for arbitrary open sets: a point
has zero $H^2$ capacity in both dimensions four and five, so puncturing
a ball does not change its plate eigenvalue.
This completes the proof of Theorem~\ref{thm:main}.

\subsection{Quantitative stability and its optimal power}\label{sec:quantitative}

The ball identity
$\Gamma_1(\Omega^*)=(\beta_n/q_n)^4\lamD(\Omega^*)^2$ and
Theorem~\ref{thm:main} give
\[
 1+\delta_\Gamma\ge(1+\delta_\lambda)^2.
\]
The sharp quantitative Faber-Krahn theorem \cite{BDPV} supplies
$\kappa_n>0$ with $\delta_\lambda\ge\kappa_n\mathcal A(\Omega)^2$,
proving \eqref{eq:quantitative}. For an unbounded finite-volume open set,
apply the bounded result to an increasing exhaustion by bounded open
subsets. Their measures and first Dirichlet eigenvalues converge, and
their normalized Fraenkel asymmetries converge as well. The latter follows
from the triangle inequality for symmetric differences and convergence of
the comparison radii.

To prove optimality, let $B$ be the unit ball and set
\[
 T_\varepsilon=\operatorname{diag}(e^\varepsilon,e^{-\varepsilon},1,\ldots,1),
 \quad E_\varepsilon=T_\varepsilon B.
\]
These ellipsoids have volume $|B|$. If $v$ is the radial first plate
eigenfunction on $B$, then $v(T_\varepsilon^{-1}x)$ is admissible on
$E_\varepsilon$. Put $M=T_\varepsilon^{-1}T_\varepsilon^{-T}$.
The Fourier transform of the zero extension of $v$ is radial.
Plancherel's theorem therefore gives
\begin{equation}\label{eq:ellipsoid}
 \frac{\Gamma_1(E_\varepsilon)}{\Gamma_1(B)}
 \le\frac{(\tr M)^2+2\tr(M^2)}{n(n+2)}
 =1+\frac{8(n+4)}{n(n+2)}\varepsilon^2+O(\varepsilon^4).
\end{equation}
With $\mathbb S^{n-1}=\partial B_1$ and $d\theta$ its surface
measure, the spherical average used here is
\[
 \frac1{n\omega_n}\int_{\mathbb S^{n-1}}
 (\theta^TM\theta)^2\,d\theta
 =\frac{(\tr M)^2+2\tr(M^2)}{n(n+2)}.
\]

The best translated comparison ball for a centered ellipsoid is centered.
Indeed, if $K=E_\varepsilon\cap(a+B)$, then
$-K=E_\varepsilon\cap(-a+B)$ and convexity gives
$(K+(-K))/2\subset E_\varepsilon\cap B$.
Brunn--Minkowski yields $|K|\le|E_\varepsilon\cap B|$.
The radial function of $E_\varepsilon$ satisfies, uniformly on the sphere,
\[
 R_\varepsilon(\theta)
 =1+\varepsilon(\theta_1^2-\theta_2^2)+O(\varepsilon^2).
\]
Consequently,
\begin{align}
 \mathcal A(E_\varepsilon)
 &=\frac1{n\omega_n}\int_{\mathbb S^{n-1}}
       |R_\varepsilon^n-1|\,d\theta\notag\\*
 &=\frac4\pi|\varepsilon|+O(\varepsilon^2).\label{eq:asymellipsoid}
\end{align}
For the last coefficient, take $\theta$ uniformly distributed on
$\mathbb S^{n-1}$ and write its first two coordinates in polar form.
Their squared radius has mean $2/n$, and the angular mean of
$|\cos2\varphi|$ is $2/\pi$. Thus
$\mathbb E|\theta_1^2-\theta_2^2|=4/(n\pi)$, where $\mathbb E$
denotes expectation.
Equations~\eqref{eq:ellipsoid} and \eqref{eq:asymellipsoid} exclude a
uniform lower bound with any power smaller than two. This completes the
proof of Corollary~\ref{cor:quantitative}.

\begin{remark}[Higher dimensions]\label{rem:higher}
Here we use the definitions of Section~\ref{sec:radial} and the transformation
\eqref{eq:Fextension} in dimension $n$.
A necessary condition for the functional comparison
\eqref{eq:functional} is $\gamma=\lambda\ell^2-K\ge0$.
To see this, choose $0\le v\in C_c^\infty(\R^n)$ with $\int v^2=1$
and put $U_R(x)=R^{-n/2}v(x/R)$. The endpoint bounds give
\[
 \int(\Delta U_R)^2=O(R^{-4}),\quad
 \|F[U_R]-\ell U_R\|_2=O(R^{-1}),\quad
 \|\nabla F[U_R]\|_2=O(R^{-1}).
\]
Letting $R\to\infty$ in \eqref{eq:functional} would then give
$-\lambda\ell^2\le-K$.

In dimension six, the series \eqref{eq:scalarseries} give
\[
 \begin{gathered}
 5.135622301<q_6<5.135622302,\quad
 5.905678235<\beta_6<5.905678236,\\
 -.01640<\gamma_6<-.01639.
 \end{gathered}
\]
These enclosures use indices $0,\ldots,23$. On $0\le r\le6$, the
absolute remainders for $j,i$ and for $b_6,b_{6,I}$ are bounded by
\[
 \frac{2\,9^{24}}{24!(3)_{24}},\quad
 \frac{2\,9^{24}}{6\cdot24!(4)_{24}},
\]
respectively. To identify the first roots, put $v=1-r^2/36$.
The alternating sum through index five gives
\[
 6b_6(r)\ge
 \frac{1319+27045v+18630v^2+36045v^3+6561v^5}{89600}>0.
\]
Thus $j$ is strictly decreasing on $(0,6]$. The cross-product identity
in Section~\ref{sec:scalar} then identifies $\beta_6$ as the first
positive root. Since $\gamma_6<0$, the functional comparison fails
in dimension six. This obstruction concerns the transformation and
leaves the higher-dimensional conjecture undecided.
\end{remark}

\section*{Statements and declarations}



\subsection*{Use of generative AI}
During the preparation of this work, the author used OpenAI ChatGPT
and Anthropic Claude to assist with developing and checking mathematical arguments and
calculations, and with editing the exposition. The author
reviewed and edited the output and takes responsibility for the final
manuscript.

\subsection*{Funding.}
This research is supported by Nazarbayev University under FDCRGP
grant No.~110326FD3205.

\end{document}